\documentclass[12pt, reqno]{amsart}
\usepackage{amsmath, amsthm, amscd, amsfonts, amssymb, graphicx, color}
\usepackage[bookmarksnumbered, colorlinks, plainpages]{hyperref}
\hypersetup{colorlinks=true,linkcolor=red, anchorcolor=green, citecolor=cyan, urlcolor=red, filecolor=magenta, pdftoolbar=true}

\newtheorem{theorem}{Theorem}[section]
\newtheorem{lemma}[theorem]{Lemma}

\newtheorem{proposition}[theorem]{Proposition}
\newtheorem{corollary}[theorem]{Corollary}
\theoremstyle{definition}
\newtheorem{definition}[theorem]{Definition}
\newtheorem{example}[theorem]{Example}

\theoremstyle{remark}
\newtheorem{remark}[theorem]{Remark}
\numberwithin{equation}{section}

\begin{document}

\setcounter{page}{1}

\title[ \MakeLowercase{$c$}-numerical range of operator products on ${\mathcal B}(H)$  ]{ \MakeLowercase{$c$}-numerical range of operator products on ${\mathcal B}(H) $ }

\author[Yanfang Zhang, \MakeLowercase{and} Xiaochun Fang]{Yanfang Zhang,$^1$ \MakeLowercase{and} Xiaochun Fang$^{2*}$}

\address{$^{1}$School of Mathematical Sciences, Tongji University, Shanghai 200092, P.R.
China.}
\email{\textcolor[rgb]{0.00,0.00,0.84}{sxzyf1012@sina.com}}

\address{$^{2}$School of Mathematical Sciences, Tongji University, Shanghai 200092, P.R.
China.}
\email{\textcolor[rgb]{0.00,0.00,0.84}{xfang@tongji.edu.cn}}




\begin{abstract}
Let $\mathcal H$ be a complex Hilbert space of dimension $\geq 2$ and $\mathfrak{B}(\mathcal H)$ be the algebra of all bounded linear operators on
$\mathcal H$.  We give the form of  surjective maps on $\mathfrak{B}(\mathcal H)$  preserving $c$-numerical range of operator products when the maps satisfy  preserving  weak zero products. As a result, we obtain the characterization of surjective maps on $M_n(\mathbb C)$  preserving $c$-numerical range of operator products. The proof of the results depends on some propositions of operators in $\mathfrak{B}(\mathcal H)$, which are of different interest.
\end{abstract} \maketitle

\section{\textbf{Introduction and preliminaries}}

Let $\mathcal H$ be a complex Hilbert space with inner product  $\langle\cdot,\cdot\rangle$, $\mathfrak{B}(\mathcal H)$ and $\mathfrak{B}_s(\mathcal H)$ be  the algebra of all bounded linear operators and  the real linear space of all self-adjoint operators on $\mathcal H$, respectively.  For $c=(c_1 ,...,c_k)\in \mathbb R^k$ and $k\leq \rm{dim}~ \mathcal H$, the $c$-numerical range and $c$-numerical radius of $A\in {\mathfrak{B}}(\mathcal H)$ are respectively defined as
 $$ W_c(A)=\left\{ \sum\limits_{j=1}^{k} c_j\langle Ae_j, e_j\rangle: \{ e_1,...,e_k\}~ \rm is~ \rm an ~\rm orthonormal ~\rm subset~ in ~\mathcal H\right\},$$
$$ r_c(A)=sup\{ |\lambda|: \lambda \in W_c(A)  \}.$$

Obviously, if $k=1$, we get the classical numerical range $W(A)$ and the numerical radius $w(A)$ of $A$. If $(c_1,...,c_k)=(1,...,1)$,  $W_c(A)$ and $r_c({A})$ reduce to $k$-numerical range and $k$-numerical radius of $A$, respectively ( \cite{CL,GR,H}). Numerical range (radius) and $c$-numerical  range (radius) are important concepts and have many applications in pure and applied mathematics, especially in quantum control and quantum information. Readers can refer to \cite{GR,SDH} for more information.

It is always of interest to characterize mappings with some special properties such as leaving certain functions, subsets or relations invariant, which are called preserving problems. There  is a considerable interest about preservers of numerical range and generalized numerical range (radius)\cite{C,C1,LT1,L}.  Li  gave the form of linear map on $ M_n(\mathbb C)$ preserving $c$-numerical range (resp. radius) \cite{LT1}. It is quite different from  linear maps preserving numerical range (resp. radius)on $\mathfrak{B}(\mathcal H)$(\cite{BH}). Hou and Di \cite{HD} gave the characterization of surjective maps  preserving numerical range of operator products on  $\mathfrak{B}(\mathcal H)~\text{and}~ \mathfrak{B}_s(\mathcal H)$ and  showed that such maps have the form $A\rightarrow \pm UAU^*$ for some unitary $U\in \mathfrak{B}(\mathcal H)$.  It is natural to ask what is the form of  surjective maps preserving $c$-numerical range of  operator products.

Let $\mathcal A,\mathcal B$ be two algebras and $\Phi: \mathcal A\rightarrow\mathcal B$ is a map. For $A,B\in \mathcal A$, $AB=0\Rightarrow\Phi(A)\Phi(B)=0$ ($AB=0\Leftrightarrow\Phi(A)\Phi(B)=0$), $\Phi$ is called a map preserving zero product (preserving zero products on both sides). The characterization of surjective maps preserving zero products on both sides on standard operator algebras was given in  \cite[Lemma 2.2]{CH}.  If $c$-numerical radius is a norm,  maps preserving  $c$-numerical radius of operator products  on $\mathfrak{B}(\mathcal H)$ are maps
preserving zero products on both sides. Then by  \cite[Theorem 2.3]{CH}, the form of maps preserving $c$-numerical radius of operator products for  $\dim(\mathcal H)\geq 3$ can get when  $r_c(\cdot)$ is a norm. The same result can also be obtained from  \cite[Theorem 3.2 \text{and} Theorem 3.3 ]{BB}.  From these, it is   possible to get the form  of  surjective maps preserving $c$-numerical range of  operator products on $\mathfrak{B}(\mathcal H)$ only when  $r_c(\cdot)$ is a norm and the dimension of $\mathcal H$ is no less than 3.

  In this article, we give the characterization of  surjective maps on $M_n(\mathbb C)$  preserving $c$-numerical  range of matrix products. Moreover, we give the definition of maps preserving weak zero products (on both sides). Then  we obtain the form of surjective maps  on $\mathfrak{B}(\mathcal H)$ preserving $c$-numerical range of operator products when the maps satisfy  preserving  weak zero products. From this result, the characterization  of surjective maps on $\mathfrak{B}(\mathcal H)$  preserving $c$-numerical range of operator products can  get easily when $r_c(\cdot)$ is a norm. Here $\dim\mathcal H$ just needs to be greater than $1$. Then the  characterization of surjective maps preserving numerical range of operator products in \cite{HD} is just a special case of this result. It is worth noting that we calculate the $c$-numerical range  of some special operators in $\mathfrak{B}(\mathcal H)$  and show a type of operators with symmetric $c$-numerical range whenever $\mathcal H$ is of finite or infinite dimension. They are of different interest.

     Next we  introduce  some notations and state the assumptions in the paper.

 (1) Assume $c=(c_1,...,c_k)\in \mathbb R^k$ and $c_1\geq...\geq c_k$ with not all $c_i$'s equal.

(2)To avoid trivial considerations, we always assume $k\geq 2$.

(3) If $\dim \mathcal H=n<\infty$, we regard $\mathfrak{B}(\mathcal H)$ as $M_n(\mathbb C)$, which is the set of all complex $n\times n$ matrices. Also we  assume $n\geq k$.

(4) $\mathcal F_1(\mathcal H)$, $\mathcal S$ and $\mathcal T$  are respectively the set of rank-1 operators, the set of operators with symmetrical $c$-numerical range and the set satisfying $\{A\in \mathcal{S}:AB\in \mathcal{S} ~{ and}~BA\in \mathcal{S}~{ for~ }~B\in \mathfrak{B}(\mathcal H)  \}$ in $\mathfrak{B}(\mathcal H)$.

The following  are our main results.
\begin{theorem}\label{3.1}
Let $M_n(\mathbb C)$ be the set of all $n\times n$ complex matrices and $\Phi:M_n(\mathbb C)\rightarrow M_n(\mathbb C)$ be a surjective map. Then  $W_c(\Phi(A)\Phi(B))=W_c(AB)$ for all $A,B\in M_n(\mathbb C) $ if and only if
  $\Phi$ has one of the following forms:

   {\rm\bf{(\romannumeral 1)}} when $c_1+c_{k}\not=0$, then there is a unitary matrix $U\in M_n(\mathbb C)$ such that  $\Phi(A)=\pm UAU^*$ for all  $A\in M_n(\mathbb C) $;

  {\rm\bf{(\romannumeral 2)}} when  there exists an  integer $p$ $(1<p<k)$ satisfying $c_p+c_{k+1-p}\not=0$ and $c_j+c_{k+1-j}=0$ for $ j=1,...,p-1$, then there are  a unitary matrix $U\in M_n(\mathbb C)$, functions $\phi: \mathcal{T}\rightarrow\{-1,1\}$ and $\phi':M_n(\mathbb C)\rightarrow\{-1,1\}$ satisfying $\phi'(A)\phi'(B)={\rm sgn}(AB)$ for any $A,B\in M_n(\mathbb C) $ where
    \[{\rm sgn}(AB)= \left\{
                     \begin{array}{cc}
                      \{1\}&\text {if}~AB\not\in\mathcal{S}\\
                        \{-1,1\}&\text {if}~AB\in\mathcal{S}  \\
                     \end{array}
                     \right. \]

   such that
     \[\Phi(A)= \left\{
                     \begin{array}{cc}
                     \phi(A)UAU^*&\text {if}~A\in \mathcal{T}\\
                       \phi'(A) UAU^*&\text {if}~A\in \mathcal{S}\setminus \mathcal{T} \\
                    \pm UAU^*&\text {if}~A\in  M_n(\mathbb C) \setminus \mathcal{S} \\
                     \end{array}
                     \right. \]
 for all  $A\in  M_n(\mathbb C) $.

      {\rm\bf{(\romannumeral 3)}} when $c_j+c_{k+1-j}=0$ for $j=1,...,k$,  then there are  a unitary matrix $U\in M_n(\mathbb C)$ and   a function $g: M_n(\mathbb C)\rightarrow \{1,-1\} $  such that
       $$\Phi(A)=g(A)UAU^* \,\,\,\,\,\,\,\,\,\,\,\, \text{for all } A\in  M_n(\mathbb C) $$
       or $$\Phi(A)=ig(A) UAU^*\,\,\,\,\,\,\,\,\,\,\,\,\text{for all } A\in  M_n(\mathbb C).$$
\end{theorem}
Next result is about  surjective maps  on $\mathfrak{B}(\mathcal H)$ preserving $c$-numerical range of operator products when the maps satisfy  preserving  weak zero products.
\begin{theorem}\label{3.2}
 Let $\mathfrak{B}(\mathcal H)$ be   the algebra of all bounded linear operators on complex Hilbert space $\mathcal H$. Suppose $\Phi:\mathfrak{B}(\mathcal H)\rightarrow\mathfrak{B}(\mathcal H)$ is a surjective map preserving weak zero products on both sides. Then $W_c(\Phi(A)\Phi(B))=W_c(AB)$ for all $A,B\in \mathfrak{B}(\mathcal H) $ if and only if   $\Phi$ has one of the following forms:

   {\rm\bf{(\romannumeral 1)}} when $c_1+c_{k}\not=0$,  there is a unitary operator $U\in \mathfrak{B}(\mathcal H)$ such that  $\Phi(A)=\pm UAU^*$ for all  $A\in \mathfrak{B}(\mathcal H) $;

   {\rm\bf{(\romannumeral 2)}} When  there exists an  integer $p$ $(1<p<k)$ satisfying $c_p+c_{k+1-p}\not=0$ and $c_j+c_{k+1-j}=0$ for $ j=1,...,p-1$, then there are  a unitary operator $U\in \mathfrak{B}(\mathcal H)$, functions $\phi: \mathcal{T}\rightarrow\{-1,1\}$ and $\phi':\mathfrak{B}(\mathcal H)\rightarrow\{-1,1\}$ satisfying $\phi'(A)\phi'(B)={\rm sgn}(AB)$ for any $A,B\in \mathfrak{B}(\mathcal H) $ where
    \[{\rm sgn}(AB)= \left\{
                     \begin{array}{cc}
                      \{1\}&\text {if}~AB\not\in\mathcal{S}\\
                        \{-1,1\}&\text {if}~AB\in\mathcal{S}  \\
                     \end{array}
                     \right. \]

   such that
     \[\Phi(A)= \left\{
                     \begin{array}{cc}
                     \phi(A)UAU^*&\text {if}~A\in \mathcal{T}\\
                       \phi'(A) UAU^*&\text {if}~A\in \mathcal{S}\setminus \mathcal{T} \\
                    \pm UAU^*&\text {if}~A\in  \mathfrak{B}(\mathcal H) \setminus \mathcal{S} \\
                     \end{array}
                     \right. \]
 for all  $A\in  \mathfrak{B}(\mathcal H) $.

 {\rm\bf{(\romannumeral 3)}} When $c_j+c_{k+1-j}=0$ for $j=1,...,k$,  then there are  a unitary operator $U\in \mathfrak{B}(\mathcal H)$ and   a function $g: \mathfrak{B}(\mathcal H)\rightarrow \{1,-1\} $  such that
       $$\Phi(A)=g(A)UAU^* \,\,\,\,\,\,\,\,\,\,\,\, \text{for all } A\in \mathfrak{B}(\mathcal H) $$
       or $$\Phi(A)=ig(A) UAU^*\,\,\,\,\,\,\,\,\,\,\,\,\text{for all } A\in  \mathfrak{B}(\mathcal H).$$
       \end{theorem}

 If  $\sum\limits_{k=1}^n c_i\not=0$, then $r_c(\cdot)$ is a norm, the following result can be obtained immediately.

 \begin{corollary}
 Suppose $\Phi:\mathfrak{B}(\mathcal H)\rightarrow\mathfrak{B}(\mathcal H)$ is a surjective map satisfying $W_c(\Phi(A)\Phi(B))=W_c(AB)$ for all $A,B\in \mathfrak{B}(\mathcal H) $. If $\sum\limits_{k=1}^n c_i\not=0$, the form of $\Phi$  is either {\rm\bf{(\romannumeral 1)}} or {\rm\bf{(\romannumeral 2)}}  in Theorem \ref{3.2}.
 \end{corollary}

The paper is organized as follows. In Section 2, we show some  properties of  $c$-numerical range, which are true whether $\mathcal H$ is finite or not. We calculate  $c$-numerical range of rank-1, some rank-2 operators in $\mathfrak{B}(\mathcal H)$  and show a type of  operators with symmetrical $c$-numerical range. In Section 3, the proof of main results will be given.

\section{ $c$-numerical range of  operators on $\mathfrak{B}(\mathcal H)$ }

In \cite{LT1}, Li showed the properties of $c$-numerical range and $c$-numerical radius on $M_n(\mathbb C)$. First we list some of them which are also true in $\mathfrak{B}(\mathcal H)$ when $\mathcal H$ is of infinite dimension.

\begin{proposition}\label{p1}
Let $c=(c_1,...,c_k)\in \mathbb{R}^k$. For $A \in \mathfrak{B}(\mathcal H) $, then
\begin{enumerate}
\item \cite{R}The set $W_c(A)$ is convex for all $A\in \mathfrak{B}(\mathcal H).$

\item $W_c(UAU^*)= W_c(A)$ for any unitary $U\in \mathfrak{B}(\mathcal H) $.

\item $W_c(\lambda A)=\lambda W_c(A)$ for any $\lambda \in \mathbb C. $

\item $W_c(\lambda I+A)=\lambda\sum\limits_{i=1}^{k}c_i+ W_c(A)$  for any $\lambda \in \mathbb C.$

\item $W_\beta c(A)=\beta W_c(A)$  for any $\beta \in \mathbb R. $

\item Suppose ${c_i}$'s are not equal. Then $W_c(A)$ is a singleton if and only if $A$ is a scalar multiple of the identity.

\item The $c$-numerical radius $r_c(A)$ is a norm on $ B(\mathcal H)$ if and only if $\sum\limits_{i=1}^{k}c_i\not =0$ and not all ${c_i}$'s are equal.

\item Suppose ${c_i}$'s are not equal. Then $W_c(A)\subset\mathbb R$ if and only if $A$ is hermitian.
\end{enumerate}
\end{proposition}

In \cite{LTS}, Li described the $c$-numerical range for self-adjoint operators especially when $\mathcal H$ is infinite dimensional. Let $S$ be a self-adjoint operator in $ \mathfrak{B}(\mathcal H)$.
If ${\dim}\mathcal H =k$, denote by $\lambda_1(S)\geq...\geq\lambda_k(S)$ the eigenvalues of $S$. If $\mathcal H$ is infinite dimensional, define $$\lambda_m(S)=\rm sup\{\lambda_m(X^*SX):~X^*X=I_m\}.$$
Also $\lambda_m$ of an infinite dimensional operator $S$ can be defined as follows. Let
$$\sigma_e(A) = \cap\{\sigma(A + F) : F \in \mathfrak{B}(\mathcal H)\rm{~ has~ finite~ rank}\}$$
be {\it the essential spectrum} of $A\in \mathfrak{B}(\mathcal H)$, and let
$$\lambda_\infty(S)=\rm sup~ \sigma_{e}(S),$$
which also equals the supremum of the set
$$\sigma(S)\backslash\{S-\mu I \rm{~has~ a~ non-trivial~ finite~ dimensional ~null~ space}\}.$$
Then $S=\sigma(S)\cap(\lambda_\infty(S),\infty)$ has only isolated points, and we arrange
the elements in descending order, say,  $\lambda_1\geq\lambda_2\geq ...$ counting multiplicities. If $S$
is infinite, then $\lambda_j(S)=\lambda_j$ for each positive integer $j$. If $S$ has $m$ elements,
then $\lambda_j(S)=\lambda_j$ for $j = 1, . . . , m$, and $\lambda_j=\lambda_\infty$ for $j > m$. Then the following property is about $c$-numerical range of self-adjoint operators in $\mathfrak{B}(\mathcal H)$.

\begin{proposition} \label{p2}(See \cite{LT1, LTS}) Let $S\in \mathfrak{B}_s(\mathcal H) $.

\begin{enumerate}
 \item If $\dim\mathcal H=k$, $W_c(S)=\left[\sum\limits_{j=1}^k {c_j} \lambda_{k+1-j}(S),\sum\limits_{j=1}^k {c_{j }}\lambda_{j}(S)\right]$.

\item If $\mathcal H$ is of  infinite dimension, {\rm ${\textbf{cl}}(W_c(S))=\left[m_c(S),M_c(S)\right]$},
where

$${ m}_c(S)={\rm inf}\left\{-\sum\limits_{j=1}^l c_j \lambda_{j}(-S)+\sum\limits_{j=1}^{k-l}c_{k+1-j} \lambda_{j}(S):~0\leq l\leq k\right\}$$
and
$${ M}_c(S)={\rm sup}\left\{\sum\limits_{j=1}^l c_j \lambda_{j}(S)-\sum\limits_{j=1}^{k-l}c_{k+1-j} \lambda_{j}(-S):~0\leq l\leq k\right\}.$$
\end{enumerate}
\end{proposition}

Next we denote \[\bar c{_j}= \left\{
                     \begin{array}{cc}
                      c_j&\text{if}~\dim\mathcal H=k\\
                       \max\{c_j,0\} &\text{if}~ \dim \mathcal H>k \\
                     \end{array}
                     \right. \]
and
                     \[ \widetilde{c}{_j}= \left\{
                     \begin{array}{cc}
                      c_j&\text{if}~\dim\mathcal H=k\\
                       \min\{c_j,0\} &\text{if}~ \dim\mathcal H>k \\
                     \end{array}
                     \right.
                     \]
                     for $ j=1,..,k $.  This will be used to describe $c$-numerical range of some operators in $\mathfrak{B}(\mathcal H)$.

\begin{proposition}\label{p3}
Suppose that $T $ is a rank-1 operator in  $\mathfrak{B}(\mathcal H)$, then
 $W_c(T)$ is an  elliptical disk  with foci ${\bar c_1} {\rm tr}(T)$, ${ \widetilde{c}_k}{\rm tr}(T)$ and minor axis $({\bar c_1}-{\widetilde{c}_k})\sqrt{\parallel T\parallel^2-\mid {\rm tr}(T)\mid^2}$ or a line segment with end points ${\bar c_1} {\rm tr}(T)$ and ${ \widetilde{c}_k} {\rm tr}(T)$.
 \end{proposition}

\begin{proof} If $\dim \mathcal H=k$, the conclusion can be found in \cite{C1}. If $\dim \mathcal H=n>k$,   adding zero to $c$,  we get  $d=({d_1},...,{d_n})$ with ${d_1}\geq...\geq{d_n}$ and $W_c(A)=W_{d}(A)$ for all $A\in \mathfrak{B}(\mathcal H) $. Obviously,  $d_1=\max\{c_1,0\}$ and $d_n=\min\{c_k,0\}$. The conclusion also can be obtained.  Next we prove it when $\mathcal H$ is  infinite dimensional.

 Assume $T=x\otimes f$ with $\parallel x\parallel=1$, then there is unit $x'$ in $[x]^{\perp}$ and complex $a,d$  such that $f=ax+dx'$. then with respect to a suitable orthogonal basis in $\mathcal H=[x]\oplus[x']\oplus[x,x']^\perp$,
$T=\left(\begin{array}{cc}
                       a& d\\
                       0 & 0\\
                     \end{array}
                   \right)\oplus 0$ with $\parallel T\parallel=\sqrt{\mid a\mid^2+\mid d\mid^2}$ and ${\rm tr }(T)=a$.
       For $0\leq \theta<2\pi$,   The real part of $e^{i\theta}T$ is denoted by $Re(e^{i\theta}T)=\frac{1}{2}(e^{i\theta}T+e^{-i\theta}T^*)$. We calculate the non-zero eigenvalues of $Re(e^{i\theta}T)$, that is
$$\lambda_1=\frac{1}{2}Re(e^{i\theta}a)+\frac{1}{2}\sqrt{(Re(e^{i\theta}a))^2+|d|^2}>0$$
and
$$\lambda_2=\frac{1}{2}Re(e^{i\theta}a)-\frac{1}{2}\sqrt{ (Re(e^{i\theta}a))^2+|d|^2}<0.$$

From Proposition \ref{p2}(2),
we know $W_c(Re(e^{i\theta}T))=[m_{c}(Re(e^{i\theta}T),M_{c}(Re(e^{i\theta}T)]$, where
 \[M_{c}(Re(e^{i\theta}T))= \left\{
                     \begin{array}{cc}
                      c_k\lambda_2 & \text{if}~c_1\leq0\\
                       c_1\lambda_1+c_k\lambda_2 &\text{if}~ c_1>0 ~\text{and}~c_k<0 \\

                      c_1\lambda_1& \text{if}~ c_k\geq0. \\
                      \end{array}
                     \right. \]

 With the notation of $\bar c_j$ and $\widetilde{c}_j$, it can be  written $M_{c}(Re(e^{i\theta}T))=\bar c_1 \lambda_1+\widetilde{c}_k\lambda_2$. Writing $a$ as $a=|a|e^{i\phi}$, then
 $$\begin{array}{rl}M_{c}(Re(e^{i\theta}T))
                       =&\bar c_1\lambda_1+\widetilde{c}_k\lambda_2\\
                       =&\frac{1}{2}(\bar c_1+\widetilde{c}_k) {\rm Re}(e^{i\theta}a)+\frac{(\bar c_1-\widetilde{c}_k)}{2}\sqrt{({\rm Re}(e^{i\theta}a))^2+|d|^2}\\
                       =&\frac{1}{2}(\bar c_1+\widetilde{c}_k)| a|\cos(\theta+\phi)
                       +\frac{(\bar c_1-\widetilde{c}_k)}{2}\sqrt{| a|^2\cos^2(\theta+\phi)+|d|^2}\\
                      \end{array}$$

It is  shown in  (\cite{C1,C2}) that  $\Omega$ is an elliptical disc centered at $pe^{i\phi}$ with foci $(p-\beta)e^{i\phi}$ and $(p+\beta)e^{i\phi}$, and semi-minor axis $\alpha$ if and only if
\begin{eqnarray}\label{2.3.a}
M_{c}(Re(e^{i\theta}\Omega))=p\cos(\theta+\phi)+\sqrt{\beta^{2}\cos^{2}(\theta+\phi)+\alpha^2}.
\end{eqnarray}
 Then we can get $W_c(T)$ is an ellipse with ${\bar c_1}{\rm tr}(T)$ and ${ \widetilde{c}_k}{\rm tr}(T)$ as foci, and\\
 $\frac{1}{2}({\bar c_1}-{ \widetilde{c}_k})\sqrt{\parallel (T)\parallel^2-|{\rm tr}(T)|^2}$ as semi-minor axis. Specially, if $|{\rm tr}(T)|=\parallel T\parallel$, $W_c(T)$ degenerates to a line segment with end points ${\bar c_1}{\rm tr}(T)$ and ${ \widetilde{c}_k}{\rm tr}(T)$.

\end{proof}

\begin{proposition} \label{p4}Let $A \in \mathfrak{B}(\mathcal H) $ be rank-2 and there is  a suitable orthonormal basis such that $A= \left(
                     \begin{array}{cc}
                       a& d \\
                       0 & b \\
                     \end{array}
                   \right)\oplus 0$ with $a,b,d\in\mathbb{ C}$.

\begin{enumerate}
\item \label{2.4.1}If $a\geq b>0$ and $d=0$,  then $W_c(A)$ is the  line segment $[\widetilde c_ka+ \widetilde{c}_{k-1}b, \bar c_1a+ \bar c_2b]$.

\item \label{2.4.2}If $a>0>b$ and $d=0$,  then $W_c(A)$ is the line segment  $[\widetilde c_ka+ \bar{c}_{1}b, \bar c_1a+\widetilde c_kb]$.

\item \label{2.4.3}If $0>a\geq b$ and $d=0$,  then $W_c(A)$ is the line segment $[\bar c_1b+ \bar c_2a,\widetilde c_kb+ \widetilde{c}_{k-1}a]$.

\item If $a,b\in \mathbb R$ and $|d|^2>4|ab|$,  then $W_c(A)$ is the  elliptical disk with foci $\bar c_1a+\widetilde{c}_kb$ and $\bar c_1b+\widetilde{c}_k a$ and with minor axis $(\bar c_1-\widetilde{c}_k)|d|$.

\item  If $\dim \mathcal H\geq 3$ and $c$ satisfying $c_j+c_{k+1-j}=0$ for $j=1,2$, and $a,b\in \mathbb{R}^+$ with $0<|d|^2<4ab$, then $W_c(A)$ is not an elliptical disk.
\end{enumerate}
\end{proposition}
\begin{proof}

 (1) For the case $\dim\mathcal H<\infty$, the conclusion can get from Proposition \ref{p2} easily. If $\dim\mathcal H=\infty$, $A$ is self-adjoint and  $W_c(A)=[m,M]$, where $M=\max\{c_1a+c_2b,~c_1a,~0\}$ and $m=\min\{c_ka+c_{k-1}b,~c_ka,~0\}$. That is
                    \[M= \left\{
                     \begin{array}{cc}
                      c_1a+c_2b &\text{if}~c_2\geq0\\
                       c_1a &\text{if}~ c_1>0 ~\text{and}~c_2<0 \\
                    0& \text{if}~c_1\leq0 \\
                      \end{array}
                     \right. \]
                    and
\[m= \left\{
                     \begin{array}{cc}
                     c_ka+c_{k-1}b & \text{if}~c_{k-1}\leq0\\
                       c_ka&\text{if}~ c_{k-1}>0 ~\text{and}~c_k<0 \\
                      0& \text{if}~c_k\geq0.\\
                      \end{array}
                     \right. \]

There are 9 combinations of $m,~M$, five of which can be possible  with $c_1\geq c_2\geq...\geq c_k$ as follows:\\
(a) if $c_k\geq 0$, then $m=0,~ M=c_1a+c_2b $;\\
(b) if $c_2>0~ \text{and}~ c_{k-1}<0$, then $m=c_ka+c_{k-1}b$ and $M=c_1a+c_2 b$;\\
(c) if $c_1\leq 0$, then  $m=c_ka+c_{k-1}b$ and $M=0$; \\
(d) if $c_1>0~ \text{and}~ c_2<0$, then  $m=c_ka+c_{k-1}b$ and $M=c_1a$;\\
(e) if $c_{k-1}>0~ \text{and}~ c_k<0$, then  $m=c_ka$ and $M=c_1a+c_2b$.

From the cases, we get $m=\widetilde{c}_ka+\widetilde{c}_{k-1}b$ and $M=\bar c_1a+\bar c_2b$.

 (2),(3) can get as the proof of (1) similarly.

(4) Assume $a,b\in \mathbb{R}$ and $|d|^2>4|ab|$.  For $0\leq \theta<2\pi$, the real part of $e^{i\theta}A$ is
 $Re(e^{i\theta}A)=\frac{1}{2}(e^{i\theta}A+e^{-i\theta}A^*)$. We calculate the non-zero eigenvalues of $Re(e^{i\theta}A)$ denoting $\lambda_1,~\lambda_2$.  Then
\begin{equation}\label{e2.4.1}
\lambda_1=\frac{(a+b)\cos \theta+\sqrt{(a-b)\cos^2\theta+|d|^2}}{2}
\end{equation}
\begin{equation}\label{e2.4.2}
\lambda_2=\frac{(a+b)\cos \theta-\sqrt{(a-b)\cos^2\theta+|d|^2}}{2}
\end{equation}
 Obviously, $\lambda_1>0>\lambda_2$. So we have  $M_{c}(Re(e^{i\theta}A))= {\bar c_1\lambda_1+\tilde{c}_k\lambda_2}$. Then
$$M_{c}(Re(e^{i\theta}A))=\frac{(\bar c_1+\widetilde{c}_k)(a+b)\cos\theta}{2}+\frac{(\bar c_1-\widetilde{c}_k)\sqrt{(a-b)^2\cos^2\theta+|d|^2}}{2}.$$

Together with Equation \eqref{2.3.a}, it follows that $W_c(A)$ is   an  elliptical disk with foci $\bar c_1a+\widetilde{c}_kb$ and $\bar c_1b+\widetilde{c}_k a$ and minor axis $(\bar c_1-\widetilde{c}_k)|d|$.

(5) Without loss of generality, suppose $a\geq b$. For $0\leq \theta<2\pi$, the non-zero eigenvalues of  $Re(e^{i\theta}A)$ is same as that \eqref{e2.4.1} and \eqref{e2.4.2}. Besides, for $|d|^2<4ab$, there must be $\phi\in(0,\frac{\pi}{2})$ such that $\cos^2\phi=\frac{|d|^2}{4ab}$. Then

\[\left\{
                     \begin{array}{cc}
                     \lambda_1\geq\lambda_2>0 &\text{if}~\theta\in S_1\\
                       \lambda_1>0>\lambda_2 &\text{if}~ \theta\in S_2\\
                      0>\lambda_1\geq\lambda_2& \text{if}~ \theta\in S_3 .\\
                      \end{array}
                     \right. \]
where $S_1=[2\pi-\phi,2\pi)\cup[0, \phi]$, $S_2= (\phi,\pi-\phi)\cup(\pi+\phi,2\pi-\phi)$ and $S_3= [\pi-\phi,\pi+\phi]$.

We know $W(Re(e^{i\theta}A))$ is interval $[m(\theta), M(\theta)]$.  From the results in  (1), (2) and (3) above, we get
\begin{equation}
M(\theta)
                       =\left\{
                     \begin{array}{cc}
                     \bar c_1\lambda_1+\bar c_2\lambda_2 &\text{if}~\theta\in S_1\\
                       \bar c_1\lambda_1+\widetilde{c}_k\lambda_2 &\text{if}~\theta\in S_2 \\
                      \widetilde{c}_{k-1}\lambda_1+\widetilde{c}_{k}\lambda_2&\text{if}~ \theta\in S_3 .\\
                      \end{array}
                     \right.
                     \nonumber
\end{equation}

 With $c_1+c_k=0$ and $c_2+c_{k-1}=0$, we get
 \[M(\theta)
                       =\left\{
                     \begin{array}{cc}
                     c_1\lambda_1 &\text{if}~\theta\in S_1\\
                        c_1(\lambda_1-\lambda_2) &\text{if}~ \theta\in S_2\\
                     - {c}_{1}\lambda_2&\text{if}~  \theta\in S_3. \\
                      \end{array}
                     \right. \]\\
   when $k=2$ and $\dim\mathcal H\geq 3$;\\
   or
                 \[M(\theta)
                       =\left\{
                     \begin{array}{cc}
                     c_1\lambda_1 +c_2\lambda_2&\text{if}~\theta\in S_1\\
                        c_1(\lambda_1-\lambda_2) &\text{if}~\theta\in S_2\\
                    - {c}_{2}\lambda_1 - {c}_{1}\lambda_2&\text{if}~ \theta\in S_3 .\\
                      \end{array}
                     \right. \]
when $k\geq3$.

    Denote    \begin{equation}
\dot{c}
                       =\left\{
                     \begin{array}{cc}
                      \frac{c_1+c_2}{2} &\text{if}~k\geq3\\
                      \frac{c_1}{2} &\text{if}~ k=2 \\
                      \end{array}
\right.
\nonumber~~~
 \text{and~~~~}~~~~~~\ddot{c}
                       =\left\{
                     \begin{array}{cc}
                      \frac{c_1-c_2}{2} &\text{if}~k\geq3\\
                      \frac{c_1}{2} &\text{if}~ k=2. \\
                      \end{array}
\right.
\nonumber
\end{equation}

$M(\theta)$ can be written as
\begin{equation}\label{2.4.a}
 M(\theta)  =\left\{
                     \begin{array}{cc}
                    \dot{c}(a+b)\cos\theta+\ddot{c}\sqrt{(a-b)^2\cos^2\theta+|d|^2 } & \text{if}~ \theta\in S_1\\
                       (\dot{c}+\ddot{c})\sqrt{(a-b)^2\cos^2\theta+|d|^2 } &\text{if}~ \theta\in S_2 \\
                     -  \dot{c}(a+b)\cos\theta+\ddot{c}\sqrt{(a-b)^2\cos^2\theta+|d|^2 }& \text{if}~  \theta\in S_3. \\
                      \end{array}
                     \right.
 \end{equation}

    If $W_c(A)$ is an ellipse, then $M(\theta)$ has  the form as  \eqref{2.3.a}, i.e. there is $p,q,t\in \mathbb R$ and $\omega\in [0, 2\pi)$ such that
    \begin{equation}\label{2.4.b}
M(\theta)=p\cos(\theta+\omega)+\sqrt{ q^{2}\cos^{2}(\theta+\omega)+t^2}.
 \end{equation}

  From  \eqref{2.4.a}, we get $M(\theta)=M(\theta+\pi)$. Then it can be deduced $p=0$ in \eqref{2.4.b}. So  \eqref{2.4.b} can be written as
    \begin{equation}\label{2.4.c}
M(\theta)=\sqrt{q^{2}\cos^{2}(\theta+\omega)+t^2}. \tag{$2.3'$}
 \end{equation}

Next we show $\omega=0$. For $\theta\in (\phi,\pi-\phi)\cup(\pi+\phi,2\pi-\phi) $,  \eqref{2.4.a} and  \eqref{2.4.c} tell that
  \begin{equation}\label{2.4.d}
q^{2}\cos^{2}(\theta+\omega)+t^2=(\dot{c}+\ddot{c})^2(a-b)^2\cos^2\theta+(\dot{c}+\ddot{c})^2 |d|^2.
 \end{equation}

 Taking the derivative about $\theta$ on both sides of  \eqref{2.4.d}, we get
  \begin{equation}\label{2.4.e}
q^{2}\sin(2\theta+2\omega)=(\dot{c}+\ddot{c})^2(a-b)^2\sin2\theta.
 \end{equation}

 Let $\theta=\frac{\pi}{2}$,  \eqref{2.4.e} shows $\sin 2\omega=0$ i.e. $\omega=0~\text{or}~\frac{\pi}{2}~\text{or}~\frac{3\pi}{2}$.
 Again take the derivative on both sides of  \eqref{2.4.d} and let $\theta=\frac{\pi}{2}$, we get $q^{2}\cos2\omega=(\dot{c}+\ddot{c})^2(a-b)^2$. So $\cos2\omega$ has to be non-negative and then $\omega=0$. Then  \eqref{2.4.c} can be written as
  \begin{equation}\label{2.4.f}
  M(\theta)=(\dot{c}+\ddot{c})\sqrt{(a-b)^{2}\cos^{2}\theta+|d|^2}. \tag{$2.3''$}
 \end{equation}

Let $\theta=0$ in  \eqref{2.4.a} and  \eqref{2.4.f}, we  have
 $$\dot{c}(a+b)+\ddot{c}\sqrt{(a-b)^2+|d|^2 }=(\dot{c}+\ddot{c})\sqrt{(a-b)^{2}+|d|^2}.$$
  which get $|d|^2=4ab$, a contradiction. So $W_c(A)$ is not an ellipse.
\end{proof}

\begin{proposition} Suppose $T,S~\in \mathfrak{B}(\mathcal H) $  are rank-1 with  $W_c(T)=W_c(S)$,  then

\begin{enumerate}
\item If $c_1+c_{k}\not=0$,  ${\rm tr}(T)={\rm tr}(S)$;
\item If $c_1+c_{k}=0$,  ${\rm tr}(T)={\rm tr}(S)$ or ${\rm tr}(T)=-{\rm tr}(S)$.
\end{enumerate}
\end{proposition}
\begin{proof}
By Proposition \ref{p3}, $W_c(T)=W_c(S)$ implies one of the following equations true:
\begin{equation}\label{2.5.a}
\left
\{
\begin{array}{lr}
\bar c_1{\rm tr}(T)=\bar c_1{\rm tr}(S)&\\
\widetilde{c}_k {\rm tr}(T)=\widetilde{c}_k {\rm tr}(S)&
\end {array}
\right.
\end{equation}
\begin{equation}\label{2.5.b}
\left
\{
\begin{array}{lr}
\bar c_1{\rm tr}(T)=\widetilde{c}_k{\rm tr}(S)&\\
\widetilde{c}_k {\rm tr}(T)=\bar c_1{\rm tr}(S)&
\end {array}
\right.
\end{equation}
With the fact that $c_1$ and $c_k$ are not 0 at the same time, we can get ${\rm tr}(T)={\rm tr}(S)$ from  \eqref{2.5.a}. For \eqref{2.5.b},  it shows ${\rm tr}(T)={\rm tr}(S)=0$ when $c_1+c_k\not=0$ and ${\rm tr}(T)=-{\rm tr}(S)$ when $c_1+c_k=0$. Hence, when $c_1+c_k\not=0$, ${\rm tr}(T)={\rm tr}(S)$. When $c_1+c_k=0$, either ${\rm tr}(T)={\rm tr}(S)$ or ${\rm tr}(T)=-{\rm tr}(S)$ holds.
\end{proof}
Next we  consider the symmetry of $c$-numerical range for    operators   in $\mathfrak{B}(\mathcal H)$.

\begin{proposition}\label{p6}
Let $A \in \mathfrak{B}(\mathcal H)$ be a self-adjoint operator  with finite ranks.  For $c=(c_1,...,c_k)$, there is an integer $p~(1<p<k)$ such that $c_p+c_{k+1-p}\not=0$ and $c_j+c_{k+1-j}=0$ for $ j=1,..., p-1$.
\begin{enumerate}
\item If {\rm rank~}$A<p$, then  $W_c(A)=-W_c(A)$;
\item  If {\rm rank~}$A=p$ and $ A$ or $-A$ is  positive, then $W_c(A)\not=-W_c(A)$.
\end{enumerate}
\end{proposition}
\begin{proof}
Assume  $A\in \mathfrak{B}(\mathcal H) $ is self-adjoint  with {\rm rank~}$A=a ~(a\leq p)$. Let $\lambda_1\geq...\geq \lambda_l$ be the positive eigenvalues  of $A$ and $\lambda_k\leq ...\leq \lambda_{k+1-(a-l)}$ be the negative ones. From Proposition \ref{p3}, there is $\xi,\eta\in \mathbb R$ such that
$W_c{(A)}= [\xi,\eta]$.

{\rm(1)} If $1\leq a<p$,  we can calculate  $\xi=\sum\limits_{j=1}^l c_{k+1-j}\lambda_{j}+\sum\limits_{j=1}^{a-l} c_{j}\lambda_{a+1-j}$
  and $\eta= \sum\limits_{j=1}^l c_{j}\lambda_{j}+\sum\limits_{j=1}^{a-l} c_{k+1-j}\lambda_{a+1-j}$.
 Obviously, $\xi+\eta=\sum\limits_{j=1}^l (c_{j}+c_{k+1-j})\lambda_{j}+\sum\limits_{j=1}^{a-l} ( c_j+ c_{k+1-j})\lambda_{a+1-j}=0$. So $W_c(A)=-W_c(A)$.

{\rm(2)} If $a= p$ and $A$ is positive, the eigenvalues of $A$, saying $\lambda_1,...,\lambda_p$, are all positive.  When  $\dim\mathcal H =k$, then $\xi=\sum\limits_{j=1}^p c_{k+1-j}\lambda_j$ and $\eta=\sum\limits_{j=1}^p c_{j}\lambda_j$. Easily, $\xi+\eta=\sum\limits_{j=1}^p c_{k+1-j}\lambda_j+ \sum\limits_{j=1}^p c_{j}\lambda_j\not=0$.

When $\mathcal H$ is infinite dimensional, we know
  \[\xi=\left\{ \begin{array}{cc}
                       \sum\limits_{j=1}^p c_{k+1-j}\lambda_{j} &\text{if}~ p\leq t~\text{and}~p<k-t \\
                      \sum\limits_{j=1}^{k-t} c_{k+1-j}\lambda_{j} & \text{if}~p\leq t~\text{and}~p>k-t \\
                      \sum\limits_{j=1}^p c_{k+1-j}\lambda_{j}&\text{if}~ p> t~\text{and}~p\leq k-t \\
                      \end{array}
                     \right. \]
and \[\eta=\left\{ \begin{array}{cc}
                       \sum\limits_{j=1}^p c_{j}\lambda_{j} &\text{if}~ p\leq t~\text{and}~p<k-t \\
                       \sum\limits_{j=1}^p c_{j}\lambda_{j} &\text{if}~ p\leq t~\text{and}~p>k-t \\
                     \sum\limits_{j=1}^t c_{j}\lambda_{j}&\text{if}~ p> t~\text{and}~p\leq k-t \\
                      \end{array}
                     \right. \]
where $t$ is the  integer satisfying $1\leq t\leq k$  and $c_{t}\geq0> c_{t+1}$. Next we show $\xi+\eta\not=0$ according to different cases of $t$:

  (a) If $p\leq t~\text{and}~p<k-t$,  then $\xi+\eta=\sum\limits_{j=1}^{p-1}(c_j+c_{k+1-j})\lambda_{j}+(c_p+c_{k+1-p})\lambda_{p}\not=0$.

  (b) If $ p\leq t~\text{and}~p>k-t$, then $\xi+\eta=\sum\limits_{j=1}^{k-t}(c_j+c_{k+1-j})\lambda_{j}+(c_{k-t+1}\lambda_{k-t+1}+...+c_{p}\lambda_{p})=c_{k-t+1}+...+c_{p}$. If we assume $\xi+\eta=0$, then $c_{k-t+1}\lambda_{k-t+1}+...+c_{p}\lambda_{p}=0$, and with $c_{k-t+1}\geq...\geq c_p\geq c_t\geq 0$, we get $c_{k-t+1}=...= c_p=...=c_t=0$. Also $c_{p}+c_{k+1-p}\not=0$, so $c_{k+1-p}\not=0$. By $p>k-t$, we can get $k+1-p\leq t$. That is $c_{k+1-p}\geq c_t=0$. But from the assumption $p\leq k+1-p$, it comes $ c_{k+1-p}\leq c_p=0$. We get a contradiction and $\xi+\eta\not=0$.

  (c) If $ p> t~\text{and}~p\leq k-t$, then $\xi+\eta=\sum\limits_{j=1}^{t}(c_j+c_{k+1-j})\lambda_{j}+(c_{k+1-(t+1)}\lambda_{t+1}+...+c_{k+1-p}\lambda_{p})=c_{k+1-(t+1)}\lambda_{t+1}+...+c_{k+1-p}\lambda_{p}<0$, so $\xi+\eta\not=0$.
This means $W_c(A)\not=-W_c(A)$ if $a=p$. As to $-A$ is positive, similar results can get.
\end{proof}

For the self-adjoint operators with rank more than $p$, the symmetry of  $c$-numerical range is associated with $c_j$ and $c_{k+1-j}$ for $j>p$, hence we failed to give a fixed conclusion.

\section{Maps preserve $c$-numerical range of  operator products on $\mathfrak{B}(\mathcal H)$ }

  In this section, first we give the definition of maps preserving weak zero products.

\begin{definition}
Let $\mathcal A,\mathcal B$ be two subalgebras of $\mathfrak{B}(\mathcal H)$  containing $\mathcal F_1(\mathcal H)$ and $\Phi: \mathcal A\rightarrow\mathcal B$ be a map. $\Phi$ is said to {\it preserve weak zero products} if for any $A$ and rank-1 $T\in \mathcal A$, $AT=0$ implies $\Phi(A)\Phi(T)=0$. $\Phi$ is said to {\it preserve weak zero products on both sides} if it  preserves weak zero products and for any $B$ and rank-1 $S\in \mathcal B$, $BS=0$ implies $\tilde{B}\tilde{S}=0$ where $ \Phi(\tilde{B})=B$ and $\Phi(\tilde{S})=S$.
\end{definition}
 Obviously, a map preserving zero products must be one preserving weak zero products. Next we give an example to show the inverse is not true.

\begin{example}
Let $\mathcal H=\mathcal H_1\oplus\mathcal H_2$ with $\dim H_1=2$ and $ U=\left(
\begin{array}{cc}
                       0& 1 \\
                       1& 0 \\
                     \end{array}
                   \right)\oplus I$ be an unitary operator. Define $\Phi:\mathfrak{B}(\mathcal H)\rightarrow \mathfrak{B}(\mathcal H)$ by
                    \[\Phi(A)= \left\{
                     \begin{array}{cc}
                      UAU^* & if \text{ rank} A\leq 1 \\
                        AU^*& if \text{rank}~ A\geq 2  \\
                     \end{array}
                     \right. \]
\end{example}
 It is easy to check that $\Phi$ is preserving left-weak zero product. Take non-zero $B_0\in \mathfrak{B}(\mathcal H) $ and let $A=\left(\begin{array}{cc}
                       0& 1 \\
                       0& 0 \\
                     \end{array}
                   \right)\oplus 0$ and  $B=\left(\begin{array}{cc}
                       1& 0 \\
                       0& 0 \\
                     \end{array}
                   \right)\oplus B_1$,  then $AB=0$. However, $\Phi(A)\Phi(B)=U(AU^*B)U^*=U(\left(\begin{array}{cc}
                       1& 0 \\
                       0& 0 \\
                     \end{array}
                   \right)\oplus 0)U^*\not=0$.

  From Theorem \ref{3.2}, it can be see that  for maps preserving $c$-numerical range of operator products on $\mathfrak{B}(\mathcal H)$, they are equivalent between preserving weak zero products and preserving zero products.

 Next we show some lemmas which are useful to the proof of main  theorems .

 \begin{lemma}\label{3.6}
 Let $A \in \mathfrak{B}(\mathcal H)$. The following conditions are equivalent.
\begin{enumerate}
\item $A$ is a rank-1 operator.
\item For any $B\in \mathfrak{B}(\mathcal H) $ with $AB \not=0$, there is   $t_B\in \mathbb C$ such that $W_c(AB)$  is an  elliptical disk {\rm (}a line segment{\rm )} with $\bar c_1 t_B$ and $\widetilde{c}_kt_B$ as  foci {\rm (}end points{\rm )}.

\end{enumerate}
 \end{lemma}
 \begin{proof}
 Since(1)${\Rightarrow}$(2) is clear, next we show (2)${\Rightarrow}$(1). Suppose rank $A\geq 2$,  then  there are linear independent $x_1$,$ x_2\in H$ such that $Ax_1\bot Ax_2$ and $\|Ax_1\|=\|Ax_2\|=1$. We finish the proof by three cases:

 ${\rm( a)}$: If $c_1+c_k\not=0$, we let $B_1=\alpha x_1\otimes A x_1 +\beta x_2\otimes A x_2$, where $\alpha>0>\beta$.  $AB_1$ is a self-adjoint operator in $\mathfrak{B}(\mathcal H)$. By proposition \ref{p4}(2), $W_c(AB_1)=[\bar c_1\beta+\widetilde{c}_k\alpha, \bar c_1\alpha+\widetilde{c}_k\beta]$.
 On the other hand, there is a complex $t_{B_1}$ such that  $W_c(AB_1)$ is a line segment with end points $\bar c_1t_{B_1}$ and $\widetilde{c}_kt_{B_1}$.  So we get
 \begin{equation}
\left
\{
\begin{array}{lr}
\bar c_1\beta+\widetilde{c}_k\alpha=\bar c_1t_{B_1} &\\
\bar c_1\alpha+\widetilde{c}_k\beta =\widetilde{c}_kt_{B_1}&
\end {array}
\right.
\nonumber
{\text or\,\,\,\,\,\,\,\,}
\left
\{
\begin{array}{lr}
\bar c_1\beta+\widetilde{c}_k\alpha=\widetilde{c}_kt_{B_1} &\\
\bar c_1\alpha+\widetilde{c}_k\beta =\bar c_1t_{B_1}&
\end {array}
\right.
\nonumber
\end{equation}
It deduces  $\bar c_1+\widetilde{c}_k=0$, impossible.

 ${\rm( b)}$: If $c_1+c_k=0$ and $c_2+c_{k-1}\not=0$, we let $B_2=\xi x_1\otimes A x_1 +\eta x_2\otimes A x_2$, where $\xi>\eta>0$ .  $AB$ is a self-adjoint operator in $\mathfrak{B}(\mathcal H)$. By proposition \ref{p4}(1), $W_c(AB_2)=[\widetilde{c}_k\xi+\widetilde{c}_{k-1}\eta, \bar c_1\xi+\bar c_2\eta]$.
  Also  there is a complex $t_{B_2}$ such that $W_c(AB_2)$ is a line segment with end points $\bar c_1t_{B_2}$ and $\widetilde{c}_kt_{B_2}$. So
 \begin{equation}
\left
\{
\begin{array}{lr}
\widetilde{c}_k\xi+\widetilde{c}_{k-1}\eta=\bar c_1t_{B_2} &\\
\bar c_1\xi+\bar c_2\eta=\widetilde{c}_kt_{B_2}&
\end {array}
\right.
\nonumber
{\text or\,\,\,\,\,\,\,\,}
\left
\{
\begin{array}{lr}
\widetilde{c}_k\xi+\widetilde{c}_{k-1}\eta=\widetilde{c}_kt_{B_2} &\\
\bar c_1\xi+\bar c_2\eta=\bar c_1t_{B_2}&
\end {array}
\right.
\nonumber
\end{equation}

It deduces  $\bar c_2+\widetilde{c}_{k-1}=0$, impossible.

 ${\rm(c)}$: If $c_1+c_k=0$ and $c_2+c_{k-1}=0$, we  let $B_3= x_1\otimes A x_1 + x_2\otimes A x_2$ when $\dim \mathcal H=2$. Then $W_c(AB_3)$ is a singleton, a contradiction. When $\dim \mathcal H\geq 3$, let $B_4= ax_1\otimes A x_1 + \varepsilon x_2\otimes A x_1+ b x_2\otimes A x_2$ where $a>b>0$ and $0<\varepsilon^2<4ab$, then  we know $W_c(AB_4)$ is not an ellipse from proposition \label{p4}(5),   a contradiction.  So $A$ must be rank-1.
 \end{proof}

  \begin{lemma} \label{3.7}
  Let $T \in \mathfrak{B}(\mathcal H) $, there is a  function $h:\mathcal H\times\mathcal H\rightarrow\{1,-1\}$ satisfying $\langle Tx,f\rangle=h(x,f)\langle x,f\rangle$ for each $x,f\in \mathcal H $, then $T=\pm I$.
  \end{lemma}
  \begin{proof}
  Firstly,  we will show  $\langle Tx,x\rangle=\langle x,x\rangle$ for any unit $x\in \mathcal H$  or $\langle Tx,x\rangle=-\langle x,x\rangle$ for any unit $x\in \mathcal H$. If not, there are   linearly independent unit $x_1,x_2\in \mathcal H$ such that  $\langle Tx_1,x_1\rangle=\langle x_1,x_1\rangle$ and  $\langle Tx_2,x_2\rangle=-\langle x_2,x_2\rangle$, then
$$\begin{array}{rl}\langle T(x_1+x_2),x_1+x_2\rangle =&\langle T x_1,x_1\rangle+\langle Tx_2,x_2\rangle+\langle Tx_1,x_2\rangle+\langle Tx_2,x_2\rangle\\
=&\langle x_1,x_1\rangle+\langle x_1,x_2\rangle-\langle x_2,x_2\rangle-\langle x_2,x_1\rangle.\\
\end{array}$$
  We know  either $\langle T(x_1+x_2),x_1+x_2\rangle= \langle x_1+x_2, x_1+x_2\rangle$ or  $\langle T(x_1+x_2),x_1+x_2\rangle= -\langle x_1+x_2, x_1+x_2\rangle$. Hence, we get $\langle x_1,x_2\rangle=\|x_1\|^2$ or $\langle x_2,x_1\rangle=\|x_2\|^2$, a contradiction with the linear independence  of $x_1,x_2$.

  Assume $\langle Tx,x\rangle=\langle x,x\rangle$  for $x\in \mathcal H$. For each $f\in \mathcal H$, there is $\dot{x}\in [x]$ and $\dot{f}\in [x]^\bot$ such that $f=\dot{x}+\dot{f}$, where $[x] $ is the subspace spanned by $x$ and $[x]^\bot$ is orthogonal subspace of $[x]$.  Then $\langle Tx,f\rangle=\langle Tx,\dot{x}+\dot{f}\rangle=\langle x, \dot{x}\rangle$. So  $\langle Tx,f\rangle=\langle x,f\rangle$   for any $x,f\in \mathcal H$. Then $T=I$. If $\langle Tx,x\rangle=-\langle x,x\rangle$ holds for $x\in \mathcal H$, $T=-I$ can be get similarly.
  \end{proof}

  \begin{lemma}\cite[Wigner Theorem]{M}{\label{3.8}}
  Let $\mathcal H $ be a complex Hilbert space and $T:\mathcal H \rightarrow \mathcal H $  be a bijective satisfying $|\langle Tx, Ty\rangle|=|\langle x ,y\rangle|$ for any $x,y\in \mathcal H $, then there are a unitary or anti-unitary operator $U:\mathcal H\rightarrow \mathcal H $ and a  function $\theta:\mathcal H \rightarrow \mathbb{C}$ with $|\theta(x)|\equiv 1$ such that $Tx=\theta(x)Ux$ for any $x\in \mathcal H$.
  \end{lemma}
  {\bf Proof of Theorem 3.4.} The sufficiency is obvious, here we only show the necessity, i.e. $\Phi$ is  a surjective map preserving weak zero products on both sides and $W_c(\Phi(A)\Phi(B))=W_c(AB)$ for all $A,B\in \mathfrak{B}(\mathcal H) $.

   \medskip{\bf Claim 1.} {\it $\Phi$ preserving rank-1 on both sides. }

   {Assume $T\in \mathfrak{B}(\mathcal H) $ is rank-1 and $\tilde{T}$ satisfying $\Phi(\tilde{T})=T$. For any $A\in \mathfrak{B}(\mathcal H)$ with $A\Phi(T)\not=0$, there is $\tilde{A}$ such that $\Phi(\tilde{A})=A$. Then $W_c(\Phi(\tilde{A})\Phi(T))=W_c(A\Phi(T))=W_c(\tilde{A}T)$. Because  $\Phi$ preserves weak zero products on both sides, then $\tilde{A}T\not=0$. Thus $W_c(A\Phi(T))$ is an ellipse (line segment) with foci(end points) $\bar{c_1}{\rm tr}(\tilde{A}T)$ and $\widetilde{c_k}{\rm tr}(\tilde{A}T)$. By Lemma \ref{3.6}, $\Phi(T)$ is rank-1. Similarly, $\tilde{T}$ is also rank-1.}
   Next we will finish the proof by two cases.\\

{\bf Case \uppercase \expandafter {\romannumeral 1} when $c_1+c_k\not=0$}.

Proof of this case will be finished by 3 steps.

{\bf Step 1.1.} {\it $\Phi$ is injective. }

 Assume  $A_1,A_2 \in \mathfrak{B}(\mathcal H)$ satisfying $\Phi(A_1)=\Phi(A_2)$.  From $W_c(\Phi(A_i)\Phi(x\otimes f))=W_c(A_ix\otimes f)$ for $i=1,2$, we can get $W_c(A_1x\otimes f)=W_c(A_2x\otimes f)$ for any $x, f\in \mathcal H$. With $c_1+c_k\not=0$ and Proposition \ref{p5}, $\langle A_1x,f\rangle=\langle A_2x,f\rangle$ for any $x,f\in \mathcal H$. So $A_1=A_2$ and $\Phi$ is injective.

{\bf Step 1.2.} {\it $\Phi$ is linear. }

Let $A_1$ and $A_2$ in $\mathfrak{B}(\mathcal{H})$ and $T\in \mathfrak{B}(\mathcal{H})$ be rank-one, we know  ${\rm tr}(\Phi(A_1+A_2)\Phi(T))={\rm tr}((A_1+A_2)T)$. So ${\rm tr}(\Phi(A_1+A_2)\Phi(T))={\rm tr}(A_1T)+{\rm tr}(A_2T)={\rm tr}((\Phi(A_1)+\Phi(A_2))\Phi(T))$. By the surjection of $\Phi$, $\Phi(T)$ can run over all rank-one operators.   So $\Phi(A_1+A_2)=\Phi(A_1)+\Phi(A_2)$. Similarly, we can check $\Phi$ is homogeneous.

{\bf Step 1.3.} {\it There is a unitary operator $U\in \mathfrak{B}(\mathcal{H})$  such that $\Phi(A)=UAU^*$ for any $A\in \mathfrak{B}(\mathcal{H})$ or  $\Phi(A)=-UAU^*$ for any $A\in \mathfrak{B}(\mathcal{H})$ .}

 Now we get that  $\Phi$ is a linear bijection  preserving rank-one operators in both directions. By \cite[Lemma 1.2]{H1}, $\Phi$ has one of the following forms:

($ \romannumeral 1$) There exist bijective linear operators $U$ and $V$ such that
$\Phi(x\otimes f)=Ux\otimes Vf $ for any $\quad~x,f\in \mathcal{H}$.

($ \romannumeral 2$) There exist bijective conjugate operators $U$ and $V$ such that
$\Phi(x\otimes f)=Uf\otimes Vx$ for any $\quad~x,f\in \mathcal{H}$.

 If ($ \romannumeral 1$) holds, that is, there exist bijective linear operators $U$ and $V$ such that $\Phi(x\otimes f)=Ux\otimes Vf $ for any $\quad~x,f\in \mathcal{H}$.

By considering $W_c((x\otimes f)^2)$ and Proposition \ref{p5}2.5(1), we  get ${\rm tr}((x\otimes f)^2)={\rm tr}(\Phi(x\otimes f)^2)$. That is $\langle x,f \rangle^2=\langle Ux, Vf \rangle^2$ for any $x,f\in \mathcal H$. It follows that $U,V$ are bounded. With Lemma \ref{3.7}, $V^*U=\pm I$.

For any $A\in \mathfrak{B}(\mathcal H)$, $W_c(Ax\otimes f)=W_c(\Phi(A)\Phi(x\otimes f))$ implies $\langle V^*\Phi(A)Ux,f\rangle=\langle A x,f\rangle $.
Hence $\Phi(A)=\pm UAU^{-1}$.

For any  unit $x\in \mathcal H$, $W_c((x\otimes x)^2)=W(( Ux\otimes (U^{-1})^{*}x)^2)$ implies $Ux\otimes (U^{-1})^*x$ is self-adjoint. So there is a positive number $r$ such that $U^{-1}x=rUx$ and $\|x\|^2=r\|Ux\|^2$ for unit $x$. Thus we get $UU^*=rI$. Let $U_1=\frac{1}{\sqrt{r}}U$ and $U_1$ is unitary  such that $\Phi(A)=\pm U_1AU_1^*$ for any $A\in \mathfrak{B}(\mathcal H)$.

If ($ \romannumeral 2$) holds,  we can get there is conjugate unitary $U$ such that $\Phi(A)=\pm UA^*U^*$ for all $A\in \mathfrak{B}(\mathcal H)$ similarly as ($ \romannumeral 1$).

Pick an orthogonal basis $\{e_j|j\in \mathfrak{J}\}$ and define a conjugate unitary operator $ J:\mathcal H\rightarrow\mathcal H$ by $Jx=\sum_{j\in \mathfrak{J}}\bar{\xi_j}e_j$ if $x=\sum_{j\in \mathfrak{J}}\xi_je_j$. Clearly, $J^2=I$, $J^*=J$ and $A^*=JA^{t}J$, where $A^{t}$ is the transpose of $A$ respect to the basis $\{e_j:j\in \mathfrak{J}\}$. Let $U_1=UJ$, Then $U_1$ is a unitary operator and $\Phi(A)=\pm U_1A^tU_1^*$ for all $A\in \mathfrak{B}(\mathcal H)$.
Thus for any $A,B\in \mathfrak{B}(\mathcal H)$, we have
$W_c(AB)=W_c(\Phi(A)\Phi(B))=W_c(U_1A^tB^tU_1^*)=W_c((BA)^t)=W_c(BA)$,
which is impossible.\\

{\bf Case \uppercase \expandafter {\romannumeral 2} when $c_1+c_k=0$}.

{\bf Step 2.1.}  {\it For any $x,y\in \mathcal H$, there are $u_x,u_y$  and $\lambda_{x,y}\in \mathbb T=\{\lambda:|\lambda|=1\}$ such that $\Phi(x\otimes y)=\lambda_{x,y} u_x\otimes u_y$ and $|\langle x,y\rangle|=|\langle u_x,u_y\rangle|$. Specially, $\lambda_{x,x}\in \{-1,~1\}$ for any $x\in \mathcal H$ or $\lambda_{x,x}\in \{-i,~i\}$ for any $x\in \mathcal H$. }

For any $x\in \mathcal{H} $, assume $\Phi(x\otimes x)=u_x\otimes v_x$ with $\| u_x\|=\| x\|$. By considering $x\otimes x$, we get $\|x\|^2W_c(x\otimes x)=\langle u_x, v_x\rangle W_c(u_x\otimes v_x)$. Hence $\|x\|^4=\langle u_x,v_x\rangle^2{\text or }-\langle u_x,v_x\rangle^2$ and $\langle u_x, v_x\rangle u_x\otimes v_x$  is self-adjoint. So there exists $t\in \mathbb C$ such that $u_x=tv_x$. Then  $\| x\|^4=\langle u_x,v_x\rangle^2={\bar t}^{-2}\|u_x\|^4$ i.e. $t\in\{-1,1,-i,i\}$. Then  there is $\lambda_{x,x}\in \{-1,~1,-i,i\}$ such $\Phi(x\otimes x)=\lambda_{x,x}u_x\otimes u_x$  for all $x\in \mathcal H$.

For any $y\in \mathcal H $, next we will show there is $\lambda_{x,y}\in \mathbb T$ such that  $\Phi(x\otimes y)=\lambda_{x,y}u_x\otimes u_y$ and $|\langle u_x, u_y\rangle|=|\langle x, y\rangle|$.
Obviously, $\Phi$ preserves orthogonality of rank-1 operators in both directions. So the ranges of $\Phi(x\otimes y)$ and $\Phi(x\otimes y)^*$ contain in span$[u_x, u_y]$,  which is the linear subspace spanned by $u_x,u_y$. Assume $\Phi(x\otimes y)=v_{x,y}\otimes w_{x,y}$. From $W_c(\Phi(x\otimes y)\Phi(x\otimes x))=W_c((x\otimes y)(x\otimes x))$, we get $\lambda_{x,x}\langle u_x,w_{x,y}\rangle W_c(v_{x,y}\otimes u_x)=\langle x,y \rangle W_c(x\otimes x)$. If $x\bot y$, obviously $u_x\bot w_{x,y}$, then $w_{x,y}$ is linearly dependant with $u_y$. If $\langle x,y \rangle\not=0$,
$\frac{\lambda_{x,x}\langle u_x,w_{x,y}\rangle}{\langle x,y \rangle}v_{x,y}\otimes u_x$ is self-adjoint, then $v_{x,y}$ is also linearly dependant with $u_x$. Similarly by considering $(x\otimes y)( x\otimes x)$, we  get $w_{x,y}$ is linearly dependant with $u_y$. So there is a complex $\lambda_{x,y}$ such that $\Phi(x\otimes y)=\lambda_{x,y}u_x\otimes u_y$.

If $x\perp y$, $|\langle u_x, u_y\rangle|=|\langle x, y\rangle|=0$. For $W_c(u_x\otimes u_y)$ is an circle disc with diameter $(\bar{c_1}-\widetilde{c_k})\|u_x\|\|u_y\|$, we  get  $|\lambda_{x,y}|=1$ from $W_c(\Phi(x\otimes y)\Phi(x\otimes x))=W_c((x\otimes y)(x\otimes x))$. If $x\not\perp y$,  considering $\Phi(x\otimes y)\Phi(x\otimes x)$, we know $\lambda_{x,y}\lambda_{x,x}\langle u_x,u_y\rangle= \langle x,y\rangle{\text or}-\langle x,y\rangle$ and $\sqrt{\|x\|^2\|y\|^2-|\langle x,y \rangle|^2}=\sqrt{|\lambda_{x,y}|^2\|u_x\|^2\|u_y\|^2-|\lambda_{x,y}\langle u_x,u_y\rangle|^2}$. Hence $|\lambda_{x,y}|=1$ and $|\langle u_x,u_y\rangle| =|\langle x,y\rangle|$.

In fact, between $\lambda_{x,x}\in \{-1,1\}  $ and $\lambda_{x,x}\in \{i,-i\}$, only one holds for any $x\in \mathcal{H}$. If not, assume there are $x_0$ and $y_0$ satisfying $\Phi(x_0\otimes x_0)= i u_{x_0}\otimes u_{x_{0}}$ and  $\Phi(y_0\otimes y_0)= u_{y_{0}}\otimes u_{y_0}$. Take $z_0\in\text{ span}[x_0,y_0]$. Considering $W_c(\Phi(x_0\otimes x_0)\Phi(z_0\otimes z_0))$, we get $\lambda_{z_0}\in\{i,-i\}$. However, we get $\lambda_{z_0}\in\{1,-1\}$ by  considering $W_c(\Phi(y_0\otimes y_0)\Phi(z_0\otimes z_0))$.  That is a contradiction.

{\bf Step 2.2.} {\it There are a unitary operator $U: \mathcal{H}\rightarrow\mathcal{H}$ and a function $g: \mathfrak{B}(\mathcal H)\rightarrow \{-1,1\}$ such that $\Phi(A)= g(A)UAU^* $ for any $A\in \mathfrak{B}(\mathcal H)$ or $\Phi(A)=ig(A)UAU^* $ for any $A\in \mathfrak{B}(\mathcal H)$.}



By Step {2.1}, there are $u_x$ with $\|u_x\|=\|x\|$  such that $\Phi(x\otimes x)=\lambda_{x,x}u_x\otimes u_x$ for any $x\in \mathcal H$.  Then define a map $V: \mathcal H\rightarrow \mathcal H$ by $ Vx=u_x$ and $V(e^{i\theta}x)=e^{i\theta}Vx$ for any $\theta\in [0,2\pi)$.  Next we show $V$ is injective. Assume there is $x,y$ such that $Vx=Vy$. By considering $(x\otimes x)(y\otimes y)$, we get $\langle y,x\rangle W_c(x\otimes y)=\langle Vy, Vx\rangle W_c(Vx\otimes Vy)$. then $x=ty$ and $t^2=1$ or $t^2=-1$. Since $V(tx)=tVx$, it   deduces $Vx\not=Vy$ for $t\not=1$. So $x=y$. Because  $\Phi$ is surjective and preserves rank-1 on both sides, we get $V$ is surjective. Then there is a bijective map $V$ with $|\langle Vx, Vy\rangle|=|\langle x,y\rangle|$such that $\Phi(x\otimes y)=\lambda_{x,y}Vx\otimes Vy$ for any $x, y\in \mathcal H$ where $\lambda_{x,y}\in \mathbb T$ may be not same as above and we still use the same symbol. By Lemma \ref{3.8}, there are a unitary or anti-unitary operator $U:\mathcal H\rightarrow\mathcal H$ and $\theta:\mathcal H \rightarrow \mathbb{C}$ with $|\theta(x)|=1$ such that $Vx=\theta(x)Ux$. So
$\Phi(x\otimes y)=\lambda_{x,y}\theta(x)\overline{\theta(y)}Ux\otimes Ux\,\,\text{for any}~ x,y\in \mathcal H.$

Let $h(x,y)=\theta(x)\overline{\theta(y)}\lambda_{x,y}$, then $\Phi(x\otimes y)=h(x,y)Ux\otimes Ux\,\,\text{for any}~ x,y\in \mathcal H.$ It is obvious that $h(x,x)=\lambda_{x,x}$. Next we show $h(x,y)\in \{-1,1\}$ for any $x,y\in \mathcal H$ or $h(x,y)\in \{-i,i\}$ for any $x,y\in \mathcal H$. Take $z\in \text{span}[x,y]$ with $\langle z,x\rangle\not=0$ and $\langle z,y\rangle\not=0$. From $W_c(\Phi(z\otimes z)\Phi(x\otimes y))=W_c((z\otimes z)(x\otimes y))$,  we get $h(z,z)h(x,y)W_c(Uz\otimes Ux)=W_c(z\otimes x)$. For $ W_c(z\otimes x)$ is an ellipse but not a circle, $h(z,z)h(x,y)=1~\text{or}-1$. With the fact $h(x,x)\in \{-1,1\}$ for all $x\in \mathcal H$ or $h(x,x)\in \{-1,1\}$ for all $x\in \mathcal H$, we get  $h(x,y)\in \{-1,1\}$ for any $x,y\in \mathcal H$ or $h(x,y)\in \{-i,i\}$ for any $x,y\in \mathcal H$. Define  a function $f:\mathcal F_1(\mathcal H)\rightarrow \{-1,1\}$ with $f(x\otimes y)=h(x,y)$, then either
\begin{equation}\label{e3.1}
\Phi(x\otimes y)=f(x\otimes y)Ux\otimes Uy~\,\,\,\,\,\text{for any}~ x,y\in \mathcal H
 \end{equation}
 or
 \begin{equation}\label{e3.2}
\Phi(x\otimes y)=if(x\otimes y)Ux\otimes Uy~\,\,\,\,\,\text{for any}~ x,y\in \mathcal H.
 \end{equation}

If \eqref{e3.1} holds, for any $A$ with $\text{rank}A\geq 2$ and $x\not\in\ker A$, $W_c((\Phi (A)Ux\otimes UAx)=W_c(Ax \otimes Ax)$.  It implies $f(x\otimes Ax)\Phi (A)Ux\otimes UAx$ is self-adjoint. Then there is $t_{x,A}\in \{-1,1\}$ such that $U^*\Phi (A)Ux= t_{x,A}Ax$. In fact we can denote $t_{x,A}=t_A$ because it  is independent with $x$. For any $x,y\not\in \ker A$, $U^*\Phi (A)Ux= t_{x,A}Ax$ and $U^*\Phi (A)Uy= t_{y,A}Ay$. On the other hand, $U^*\Phi (A)U(x+y)= t_{x+y,A}A(x+y)$. So $(t_{y,A}-t_{x+y,A})Ay=(t_{x+y,A}-t_{x,A})Ax$. $t_{x,A}\not=t_{y,A}$ will deduce $Ax=0$ or $Ay=0$, impossible. Define a function $g:\mathfrak{B}(\mathcal H)\rightarrow \{-1,1\}$ with \[g(A)= \left\{
                     \begin{array}{cc}
                      f(A)&\text{rank}A=1\\
                       t_A & \text{rank}A\geq 2, \\
                     \end{array}
                     \right. \] then $\Phi (A)= g(A)U A U^*$ for any $A\in \mathfrak{B}(\mathcal H)$.

If \eqref{e3.2}  holds, $-i\Phi$ is accord with \eqref{e3.1}. Then $\Phi (A)= ig(A)U A U^*$ for any $A\in \mathfrak{B}(\mathcal H)$.

Similarly as  Step 1.3 of the proof in  {\bf Case \uppercase \expandafter {\romannumeral 1}}, it is impossible that $U$ is anti-unitary.

{\bf Step 2.3.} {\it $\Phi$ has the form in Theorem \ref{3.2}.}

When  there exists an  integer $p$ $(1<p<k)$ satisfying $c_p+c_{k+1-p}\not=0$ and $c_j+c_{k+1-j}=0$ for $ j=1,...,p-1$, the form $\Phi (A)= ig(A)U A U^*$ for $A\in \mathfrak{B}(\mathcal H)$ is impossible. If not, let $D\in  \mathfrak{B}(\mathcal H)$ is a project with rank~$D=p$.  $W_c(\Phi(D)^2)=W_c(D^2)$ deduces
$-W_c(D)=W_c(D)$, which is contradict with Proposition \ref{p6}.

Remembering  $\mathcal{S}=\{A\in \mathfrak{B}(\mathcal H):~ W_c(A)=-W_c(A)\}$ and $\mathcal{T}=\{A\in \mathcal{S}:~AB\in \mathcal{S} ~{\rm and}~BA\in \mathcal{S}~{\rm for~ any }~B\in \mathfrak{B}(\mathcal H) \}$, we know $W_c(\Phi(I)\Phi(B))=W_c(B)$ for any operator $B\not\in \mathcal{S}$. Thus $g(B)=g(I)=1 ~\text{or}-1$ for any operator $B\in \mathfrak{B}(\mathcal H)\setminus\mathcal{S}$.

  For any $T\in \mathcal{T}$, define the function  $\phi: \mathcal{T}\rightarrow\{-1,1\}$ with $\phi(T)=g(T)$. For any $S\in \mathcal{S}\setminus\mathcal{T}$, define the function  $\phi':{B}(\mathcal H)\rightarrow\{-1,1\}$ with $\phi'(S)=g(S)$. Because $W_c(\Phi(S)\Phi(A))=W_c(SA)$ for any $A\in \mathfrak{B}(\mathcal H) $, $\phi'$ must  satisfy $\phi'(S)\phi'(A)={\rm sgn}(SA)$ for any $A\in \mathfrak{B}(\mathcal H) $ where
    \[{\rm sgn}(SA)= \left\{
                     \begin{array}{cc}
                      \{1\}&SA\not\in\mathcal{S}\\
                        \{-1,1\}&SA\in\mathcal{S}.  \\
                     \end{array}
                     \right.\]

 The proof is finished.\hfill$\Box$

   \begin{remark}
   For the form  {\rm\bf{(\romannumeral 3)}} in Theorem \ref{3.2}, it is not easy to show all the elements in $\mathcal S$ and $\mathcal T$. Here we show some subsets of them.
   \end{remark}
   \begin{lemma}
  For $c=(c_1,...,c_k)$, if there exists an  integer $p$ $(1<p<k)$ satisfying $c_p+c_{k+1-p}\not=0$ and $c_j+c_{k+1-j}=0$ for $ j=1,...,p-1$. Let $V_1=\{T: T ~\text{is a positive or non-positive operator with }$ ${\rm rank} ~T= p~ \}$ and $V_2=\{T: T~\text{is a}$ $\text {self-adjoint operator with}~{\rm rank} ~T<p \}$, then $V_1\subset \mathfrak{B}(\mathcal H)\setminus \mathcal S$, $V_2\subset \mathcal S$, $\mathcal F_{1}(\mathcal H)\subset \mathcal T $ and $\mathcal T\not=\mathcal S$.
   \end{lemma}
   \begin{proof}
From Proposition \ref{p6}(1),(2), we know $V_1\subset \mathfrak{B}(\mathcal H)\setminus \mathcal S$, $V_2\subset \mathcal S$.    For  each  rank-1  operator $T$ and any operator $A\in\mathfrak{B}(\mathcal H) $, it comes  $W_c(AT)=-W_c(AT)$. Then $\mathcal F_{1}(\mathcal H)\subset \mathcal T $.
 Let $D=I_{p-1}\oplus(-1)\oplus 0$ and $G=I_{p-1}\oplus(-2)\oplus 0$.   $W_c(D)=[\sum\limits_{j=1}^{p-1}c_{k+1-j}-c_1,~\sum\limits_{j=1}^{p-1}c_{j}-c_k]$. Obviously,  $W_c(D)=-W_c(D)$. Then  $D\in \mathcal S$. However $DG$ is an positive operator with rank~$(DG)=p$, then $W_c(DG)\not=-W_c(DG)$. So $D\not\in \mathcal T$ and $\mathcal T\not=\mathcal S$.
\end{proof}

  At the end, we finish the proof of Theorem \ref{3.1}.

{\bf Proof of Theorem 1.1.}

{\bf Claim 1.} {\it $\Phi(A)=0\Leftrightarrow A=0.$}

If $\sum\limits_{j=1}^kc_j\not=0$, then $W_c(A)=\{0\}$ if and only if $A=0$ and the claim is true. If $\sum\limits_{j=1}^kc_j=0$, $W_c(T)=\{0\}$ implies $T$ is multiple of the identity. Assume $\Phi(A)=0$, then $AB$ is a multiple of the identity for any $B\in M_n(\mathbb C)$. Let $B$ run over all rank-1 matrices, we get $AD=0$ for any rank-1 matrices. Because any matrix in $M_n(\mathbb C)$ is a combination of rank-1 matrices. Thus $AB=0$ for all $B\in M_n(\mathbb C)$. Then $A=0$.

{\bf Claim 2.} {\it $\Phi$ preserves zero products on both sides.}

If $\sum\limits_{j=1}^kc_j\not=0$, obviously $\Phi$ preserves zero products. Next we show if $\sum\limits_{j=1}^kc_j=0$, $\Phi$ also preserves zero products. When $\sum\limits_{j=1}^kc_j=0$, On contrary assuming that  there are non-zero $A,B\in M_n(\mathbb C)$ with $AB=0$ such that $\Phi(A)\Phi(B)=tI\not=0$.   Then $\Phi(A)$ and $\Phi(B)$ are both invertible. Taking $x,f\in \mathbb C^n$ such that $x\in \ker A$ and  $f\not\in \ker{A}^*$, we see $Ax\otimes f=0$ and $x\otimes fA\not=0$. This deduces $\Phi(A)\Phi(x\otimes f)=sI$ where $s\in \mathbb C$. If $s=0$, we get $\Phi(x\otimes f)=0$, impossible. If $s\not=0$, then $\Phi(x\otimes f)\Phi(A)=tI$  and $x\otimes fA=0$, impossible. So $AB=0$ implies $\Phi(A)\Phi(B)=0$.

So from Theorem \ref{3.2}, $\Phi$ has the form in Theorem \ref{3.1}. \hfill$\Box$

\bibliographystyle{amsplain}

\end{document}